\documentclass[preprint,12pt]{elsarticle}
\usepackage{indentfirst}
\usepackage{natbib}
\usepackage{mathrsfs}
\usepackage{amsfonts}
\usepackage{graphicx}
\usepackage{rotating}
\usepackage{enumerate}
\usepackage{multirow}
\usepackage{amsmath}
\usepackage{amsthm}

\makeatletter
\@addtoreset{equation}{section}\makeatother

\theoremstyle{definition}
\newtheorem{Theorem}{Theorem}[section]
\newtheorem{Lemma}[Theorem]{Lemma}
\newtheorem{Corollary}[Theorem]{Corollary}
\newtheorem{Proposition}[Theorem]{Proposition}

\newtheorem{Assumption}{Assumption}[section]

\makeatletter
\@addtoreset{equation}{section}\makeatother

\begin{document}
\begin{frontmatter}
\title{Consistency of Bayesian Linear Model Selection With a Growing Number of Parameters}
\author{Zuofeng Shang and Murray K. Clayton}

\address{Department of Statistics\\
         University of Wisconsin, Madison\\
         Madison, WI 53706}

\begin{abstract}
Linear models with a growing number of parameters have been widely
used in modern statistics. One important problem about this kind of
model is the variable selection issue. Bayesian approaches, which
provide a stochastic search of informative variables, have gained
popularity. In this paper, we will study the asymptotic properties related to 
Bayesian model selection when the model dimension $p$ is growing with the sample size $n$. 
We consider $p\le n$ and provide sufficient conditions under which: (1) with large probability,
the posterior probability of the true model (from which samples are drawn) uniformly dominates the
posterior probability of any incorrect models; and (2) 
the posterior probability of the true model converges to one in probability.
Both (1) and (2) guarantee that the true model will be selected
under a Bayesian framework. We also demonstrate several situations
when (1) holds but (2) fails, which illustrates the difference between
these two properties. Finally, we generalize our results to include
$g$-priors, and provide simulation examples to 
illustrate the main results.
\end{abstract}

\end{frontmatter}

\section{Introduction}

This work was motivated by efforts to analyze remotely sensed
(satellite) data which consists of multiple spatial images. In the
setting of interest, one image corresponds to a ``response'' while
others correspond to covariates.  To find the relationship between the
response and covariate spatial images, Zhang {\it et al.} (2010) proposed a
functional concurrent linear model with varying coefficients and
applied a wavelet approach to transform this model into a linear model
(with a particular design matrix) which contains an $n$-vector of
responses and a sparse $p$-vector of wavelet coefficients. Since the
images contain thousands of pixels, the model dimension $p$, which is
determined by the maximum decomposition level in the wavelet
expansion, has to be large so that sufficiently fine details in the
target images can be captured. On the other hand, $p$ has an upper
bound $p \leq (K+1)n$, where $K$ is the total number of covariate images 
involved in the model. This is because each spatial image corresponds 
to a vector of wavelet coefficients which has dimension not exceeding $n$, 
and there are $K+1$ images in total with one of them representing the intercept 
and others the slopes.  An important question is how to select the
nonzero coefficients in the model, which is essentially a variable selection
problem. Zhang {\it et al.} (2010) adopted  a Lasso approach to address
this.

The problem they handle relies on a specific design matrix induced by
the wavelet structure.  It is of interest, to frame the variable
selection problem more broadly.  More precisely, we suppose that data
are drawn from the linear model
\begin{equation}\label{basic:linear:model}
\textbf{y}=X\beta+\epsilon,
\end{equation}
where $\epsilon\sim N(\textbf{0},\sigma_0^2 I_n)$ is an $n$-vector of
errors, $\textbf{y}=(y_1,\ldots,y_n)^T$ is an $n$-vector of responses,
$\beta=(\beta_1,\ldots,\beta_p)^T$ is a $p$-vector of parameters and
$X=(X_1,\ldots,X_p)$ is a $n\times p$ design matrix with $X_j$ the $j$th column of $X$. 
It is also assumed that only a subset of 
$X_1,\ldots,X_p$ contribute to $\textbf{y}$ and we are interested in 
selecting the variables in this subset.  
 
We consider a Bayesian variable selection (BVS) approach based on
model (\ref{basic:linear:model}). The Bayesian model to be considered
is a variation of George and McCulloch (1993) and has been studied by Clyde {\it et al.} (1998),
Clyde and George (2000), and Wolfe {\it et al.} (2004).  
Clearly, each subset of
$X_1,\ldots,X_p$ defines a candidate model, so there are $2^p$ of them
in total. According to George and McCulloch (1993), all the marginal
posterior probabilities of these $2^p$ models can be calculated and
the model with the largest posterior probability can be selected as the
``best" model. This motivates the formal definition of posterior model consistency (PMC).
We say that PMC holds if the true model, defined as the model from which
samples are drawn, has a posterior probability approaching one. 
Since the sum of the posterior probabilities of all models
equals one, when PMC holds, the posterior probability of any 
incorrect model will go to zero when $n$ goes to infinity so that the true model can
be correctly selected. 

PMC has been theoretically verified when $p$ is fixed 
(see Fern\'{a}ndez {\it et al.}, 2001; Moreno and Gir\'{o}n, 2005; 
Liang {\it et al.}, 2008; Casella {\it et al.}, 2009). 
However, fewer results have been derived when $p$ is
growing with $n$, an interesting and important scenario. 
For increasing $p$, Berger {\it et al.} (2003), Moreno {\it et al.} (2010) 
and Gir\'{o}n {\it et al.} (2010) proved consistency for Bayes factors.
Although PMC and consistency of Bayes factors are equivalent for fixed $p$
(see Liang {\it et al}, 2008; Casella {\it et al.}, 2009), they are different for growing $p$. 
Actually, we will see below that consistency of the Bayes factor 
is equivalent to consistency of the posterior odds ratio under a 
general setting, but that the latter form of consistency is weaker 
than PMC. Therefore, it seems valuable to separately study PMC.

In this paper we will consider two classes of design matrix $X$, both with
$p\le n$, although our results can be generalized to $p\gg n$
when combined with certain dimension reduction approaches.
In the first case, $X$ is quite general. A representative situation is that
the eigenvalues of $X^TX/n$ are uniformly bounded both above and below. Consistency
is examined when $p$ grows slower than $n$, say, $p\log{n}=o(n)$. We
find that the posterior odds in favor of any incorrect model uniformly
converges to zero, and the posterior probability of the true model
converges to one. A second case we consider occurs when $X^TX/n$ is
the identity matrix, i.e., $X^TX=nI_p$, and $p$ grows as fast as $n$,
say $p=n$.  In that case, consistency of the posterior odds ratio and
PMC are examined, i.e., the posterior odds ratio in favor of any incorrect 
model uniformly converges to zero, and the posterior probability of the true model
converges to one.  We also demonstrate how consistency of the posterior 
odds ratio can hold even though PMC fails. Finally, we generalize our results
to a $g$-prior setting proposed firstly by Zellner (1986). 

The remainder of this paper is organized as follows. In Section 2,
preliminaries and main results will be provided. In Section 3, a
numerical example related to the results of Section 2 is displayed. 
Section 4 contains the conclusion.
And technical arguments are included in Section 5.

\section{Preliminaries and main results}

Suppose the $n$ dimensional response vector
$\textbf{y}=(y_1,\ldots,y_n)^T$ and the $n$ by $p$ covariate matrix
$X=(X_1,\ldots,X_p)$ are linked by the model
\begin{equation}\label{model}
\textbf{y}=X\beta+\epsilon,
\end{equation}
where the $X_j$'s are $n$-vectors, $\beta=(\beta_1,\ldots,\beta_p)^T$
is an unknown $p$-vector and $\epsilon$ is a vector of random errors. 
Here, $X$ is allowed to be either (1) random but independent of $\epsilon$ 
or (2) deterministic. For $1\le j\le p$,
define the state variable of $\beta_j$ by $\gamma_j=I(\beta_j\neq 0)$
and $\gamma=(\gamma_1,\ldots,\gamma_p)^T$, where $I(\cdot)$ is the
indicator function. We call $\gamma$ the state vector of $\beta$ and  
denote the number of $1$'s in
$\gamma$ 
by $|\gamma|$.
The state vector $\gamma$ completely determines the
inclusion or exclusion of $\beta_j$'s in model (\ref{model}), and
therefore, can define a model
$\textbf{y}=X_\gamma\beta_\gamma+\epsilon$, where $X_\gamma$ is an
$n\times |\gamma|$ submatrix of $X$ whose columns are indexed by the
nonzero components of $\gamma$, and $\beta_\gamma$ is the subvector
(with size $|\gamma|$) of $\beta$ indexed by the nonzero components of
$\gamma$. It is natural, therefore, to call each $\gamma$ a model.  
Note that there are $2^p$ such $\gamma$'s representing $2^p$ different
models. For any state vectors $\gamma$ and $\gamma'$, let
$(\gamma\backslash\gamma')_j=I(\gamma_j=1,\gamma'_j=0)$ denote the
difference (which is also a state vector) between $\gamma$ and
$\gamma'$, i.e., the 0-1 vector indicating the variables that are
present in $\gamma$ but absent in $\gamma'$. We say that $\gamma$ is
nested in $\gamma'$ (denoted by $\gamma\subset\gamma'$) if
$\gamma\backslash\gamma'=0$. Denote the true model coefficient vector by
$\beta^0$ and the corresponding state vector by $\gamma^0$, and let
$s_n=|\gamma^0|$ denote the size of the true model.

In this paper we consider the following hierarchical Bayesian 
model which is a variation of the model used by George and McCulloch (1993)
\begin{eqnarray}\label{full-Bayesian-model}
&&\textbf{y}|\beta,\sigma^2\sim N(X\beta,\sigma^2 I_n),\nonumber\\
&&\beta_j|\gamma_j,\sigma^2\sim (1-\gamma_j)\delta_0+
\gamma_j N(0,c_j\sigma^2),\nonumber\\
&& 1/\sigma^2\sim \chi_\nu^2,\nonumber\\
&& \gamma\sim p(\gamma),
\end{eqnarray}
where $\delta_0$ is point mass measure concentrated at
zero. Hereafter, $\nu$ will be fixed  a priori. Let
$\Sigma=\textrm{diag}(\textbf{c})$ with $\textbf{c}=(c_j)_{1\le j\le p}$ a $p$-vector of
positive components, and let $\Sigma_\gamma$ be the $|\gamma|\times
|\gamma|$ sub-diagonal matrix of $\Sigma$ corresponding to
$\gamma$. Let $Z=(\textbf{y},X)$ denote the full data set. It follows
by integrating out $\beta$ and $\sigma$ that the posterior
distribution of $\gamma$ is given by
\begin{equation}\label{eq:post} p(\gamma|Z)\propto (2\pi)^{-n/2}\det(W_\gamma)^{-1/2}p(\gamma)
  \left\{\frac{2}{1+\textbf{y}^T(I_n-X_\gamma U_\gamma^{-1} X_{\gamma}^T)\textbf{y}}\right\}^{(n+\nu)/2},
\end{equation}
where $U_\gamma=\Sigma_\gamma^{-1}+ X_{\gamma}^TX_\gamma$ and
$W_\gamma=\Sigma_\gamma^{1/2}U_\gamma \Sigma_\gamma^{1/2}$. In
particular, if $\gamma=\emptyset$ (the null model containing no covariate variables), 
(\ref{eq:post}) still holds
if we adopt the conventions that $X_{\emptyset}=0$ and
$\Sigma_\emptyset=U_\emptyset=W_\emptyset=1$.  

Define $S_1=\{\gamma|\gamma^0\subset\gamma, \gamma\neq\gamma^0\}$ and
$S_2=\{\gamma|\gamma^0\,\, \textrm{is not nested in}\,\, \gamma\}$. It
is clear that $S(n)$ defined by $S(n)=S_1\cup S_2 \cup\{\gamma^0\}$ is
the class of all state vectors.  In particular, when $\gamma^0=\emptyset$,
$S_2$ is empty, and hence $S_1$ is the class of all state vectors
excluding $\gamma^0$. As was found by Liang {\it et al.} (2008), we will see
later in this section that whether $\gamma^0$ is null or nonnull will
result in some differences in the main results (especially in the
assumptions that are needed to establish our main results); thus, we
will treat these cases separately. When $\gamma^0$ is nonnull, we
denote $\varphi_{\min}(n)=\min\limits_{\gamma\in
  S_2}\lambda_{-}\left(\frac{1}{n}X_{\gamma^0\backslash\gamma}^T
  (I_n-P_\gamma)X_{\gamma^0\backslash\gamma}\right)$ and
$\varphi_{\max}(n)=\max\limits_{\gamma\in S_2}\lambda_{+}
\left(\frac{1}{n}X_{\gamma^0\backslash\gamma}^TX_{\gamma^0\backslash\gamma}\right)$,
where $P_\gamma=X_\gamma(X_{\gamma}^TX_\gamma)^{-1}X_{\gamma}^T$ is a
projection matrix, $\lambda_{-}(A)$ and $\lambda_{+}(A)$ are the
minimal and maximal eigenvalues of the square matrix $A$.
We also adopt the convention that $P_\emptyset=0$. For
the case that $\gamma^0=\emptyset$, both $\varphi_{\min}$ and
$\varphi_{\max}$ are meaningless, and $S_1$ will be focused on in this situation. 

Before proceeding further, we introduce several types of consistency 
central to this work. Generally speaking, to make a correct model selection
\begin{equation}\label{mpp}
\max\limits_{\gamma\neq\gamma^0}\,\,p(\gamma|Z)/p(\gamma^0|Z)\rightarrow0
\end{equation} 
should hold as $n\rightarrow\infty$, which means that the posterior probability
of the true model asymptotically dominates that of any incorrect
model. Following a framework similar to that of Zellner (1978), the term 
$p(\gamma|Z)/p(\gamma^0|Z)$, which is called the posterior odds
ratio in favor of $\gamma$, satisfies the relationship
\begin{equation}\label{BayesFactor:PMC}
p(\gamma|Z)/p(\gamma^0|Z)=BF(\gamma:\gamma^0)\frac{p(\gamma)}{p(\gamma^0)}, 
\end{equation}
where $BF(\gamma:\gamma^0):=p(Z|\gamma)/p(Z|\gamma^0)$ is the Bayes
factor of $\gamma$ versus $\gamma^0$ and $p(\gamma)/p(\gamma^0)$ is
the prior odds ratio in favor of $\gamma$. The Bayes factor is
consistent if for any $\gamma\neq\gamma^0$,
$BF(\gamma:\gamma^0)\rightarrow0$. The posterior odds ratio is
consistent if for any $\gamma\neq\gamma^0$,
$p(\gamma|Z)/p(\gamma^0|Z)\rightarrow0$. It is easy to see that
property (\ref{mpp}) implies consistency of the posterior odds ratio.
We say that posterior 
model consistency (PMC) holds if $p(\gamma^0|Z)\rightarrow1$. 
These types of consistency all have been useful in Bayesian model selection.
Representative references include (1) assessment of posterior odds ratio: 
Jeffreys (1967), Zellner (1971, 1978); (2) performance of Bayes factor: 
Berger and Pericchi (1996), Moreno {\it et al.} (1998, 2010), Casella {\it et al.} (2009); 
(3) PMC: Fern\'{a}ndez {\it et al.} (2001), Liang {\it et al.} (2008).

It is easy to see that when
\begin{equation}\label{PriorOdds}
\tilde{c}^{-1}\le\min\limits_\gamma p(\gamma)/p(\gamma^0)\le \max\limits_\gamma p(\gamma)/p(\gamma^0)\le \tilde{c}
\end{equation} 
holds for some positive constant $\tilde{c}$, consistency of the Bayes factor
is equivalent to consistency of the posterior odds ratio, and that
both are weaker than (\ref{mpp}). A special case is that $p(\gamma)=2^{-p}$ for
all $\gamma$'s, which results in an indifference prior distribution for $\gamma$, see, e.g.,
Smith and Kohn (1996).

To illustrate the relationship between PMC and (\ref{mpp}), 
note that
\begin{equation}\label{PMC:BayesFactor}
p(\gamma^0|Z)=\frac{1}{1+\sum\limits_{\gamma\neq\gamma^0}p(\gamma|Z)/p(\gamma^0|Z)},
\end{equation}
and thus $p(\gamma^0|Z)\rightarrow1$ will imply (\ref{mpp}). When $p$
is fixed, it has been noted by Liang {\it et al.} (2008) that (\ref{mpp})
implies PMC. However, when $p$ grows with $n$, it will be shown later that this may
not be true. This somewhat illustrates the difference between PMC and
(\ref{mpp}). 

In what follows, we introduce some regularity conditions that are
useful to establish our main results. We will also demonstrate some
particular situations when these conditions are satisfied.

\begin{Assumption}\label{assump:1}
  There exists a constant $C_0>0$ such that for any
  $n$, $\max\limits_{\gamma\in S(n)} p(\gamma)/p(\gamma^0)\le C_0$.
\end{Assumption}

\begin{Assumption}\label{assump:2} 
  There exist positive constants $C_1, C_2$ such that with probability
  equal to one, $\liminf\limits_n \varphi_{\min}(n)\ge C_1$ and
  $\limsup\limits_n \varphi_{\max}(n)\le C_2$.
\end{Assumption}

\begin{Assumption}\label{assump:3} There exists a positive sequence $\psi_n$ such that
  $\min\limits_{j\in\gamma^0}|\beta_j^0|\ge\psi_n$ and, as $n\rightarrow\infty$,
  $\psi_n\sqrt{n}\rightarrow\infty$.
\end{Assumption}

\begin{Assumption}\label{assump:4} 
  $p_n\rightarrow\infty$, $s_n\le p_n\le n$ and $p_n\log{n}=o(n\log(1+\min\{\psi_n^2,1\}))$.
\end{Assumption}

\begin{Assumption}\label{assump:5} 
$p_n\rightarrow\infty$, $s_n\le p_n\le n$ and $p_n\log{p_n}=o(n)$.
\end{Assumption}

Hereafter, unless otherwise explicitly stated, we will drop the subscript from $p_n$.

\begin{Assumption}\label{assump:6} There is a positive sequence
  $\bar{\phi}_n=O(n^{\delta_0})$ for some $\delta_0>0$ such that
  $\max\limits_{1\le j\le p}c_j\le \bar{\phi}_n$, where the $c_j$'s
  are the hyperparameters (in model (\ref{full-Bayesian-model}))
  controlling the prior variances of the nonzero $\beta_j$'s.
\end{Assumption}

\begin{Assumption}\label{assump:7} There is a positive sequence $\underline{\phi}_n$ such that 
  $k_n=O(\underline{\phi}_n)$ and $\min\limits_{1\le j\le p}c_j\ge
  \underline{\phi}_n$, where $k_n=\|\beta^0_{\gamma^0}\|_2^2$.
\end{Assumption}

\begin{Assumption}\label{assump:8} 
There exist $C_3>0$ and $\delta\ge 0$ such that $n^{1-\delta}\underline{\phi}_n\rightarrow\infty$, 
and for any $n$, with probability equal to one,
\begin{equation}\label{cond_viii}
\inf\limits_{\gamma\in S_1}\lambda_{-}\left(\frac{1}{n}X_{\gamma\backslash\gamma^0}'(I_n-P_{\gamma^0})X_{\gamma\backslash\gamma^0}\right)\ge C_3 n^{-\delta}.
\end{equation}
\end{Assumption}

{\bf Remark 2.1.}\,\,
\begin{enumerate}[(a).]
\item Assumption \ref{assump:1} is satisfied by some commonly used
  priors $p(\gamma)$, such as the flat prior $p(\gamma)=2^{-p}$ (Smith and Kohn, 1996). 
  More generally, if $p(\gamma_j=1)=\theta_j$
  is such that both
  $\prod\limits_{j\in\gamma\backslash\gamma^0}\left(\frac{\theta_j}{1-\theta_j}\right)$
  and $\prod\limits_{j\in
    \gamma^0\backslash\gamma}\left(\frac{1-\theta_j}{\theta_j}\right)$
  are bounded, then Assumption \ref{assump:1} is satisfied.

\item We use Assumption \ref{assump:3}  to prove
  consistency for a growing $p$. Fan and Peng (2004)
  introduced a similar assumption in the framework of 
  smoothly clipped absolute deviation (SCAD) penalized
  optimization where $\sqrt{n}$ in Assumption \ref{assump:3} was replaced by
  $1/\lambda_n$ with $\lambda_n$ the penalty parameter. This condition
  requires  the true parameters to be away from zero.
  Otherwise, it is impossible to distinguish between zero and nonzero
  parameters. 

\item Assumptions \ref{assump:4} and \ref{assump:5} define a rate on
  the dimension $p$. In particular, when $\inf\limits_n \psi_n>0$,
  Assumption \ref{assump:4} is satisfied if $s_n\le p$ and
  $p\log{n}=o(n)$. The results hold when $s_n$ is either bounded or growing with $n$.

\item Assumption \ref{assump:6} excludes the possibility
  that $\bar{\phi}_n$ is extremely large, e.g., we exclude the situation that $\bar{\phi}_n=\exp(n^\omega)$ 
  for some $\omega>0$. Assumption \ref{assump:7} requires that $\underline{\phi}_n$ is not growing slower than 
  $k_n=\|\beta^0_{\gamma^0}\|_2^2$. When the design
  matrix $X$ is nonorthogonal, we use this assumption to
  facilitate the proof of consistency (see Theorem \ref{thm1}
  below). But when $X$ is orthogonal, this assumption is
  redundant and can be removed (see Corollary \ref{cor:1} below).\qed
\end{enumerate}

Assumptions \ref{assump:1},
\ref{assump:3}--\ref{assump:7} are easily satisfied. 
The following proposition demonstrates that a
broad class of design matrices $X$ can satisfy Assumptions 2.2
and 2.8.

\begin{Proposition}\label{matrix:class}
  If the $n\times p$ matrix $X$ satisfies
  $\lambda_{-}\left(\frac{1}{n}X^TX\right)\ge c$, where $c>0$ is
  constant, then for any $\gamma\subset\bar{\gamma}$ and
  $\gamma\neq\bar{\gamma}$,
\begin{equation}\label{eq:matrix:class}
\lambda_{-}\left(\frac{1}{n}X^T_{\bar{\gamma}\backslash\gamma}(I_n-P_\gamma)X_{\bar{\gamma}\backslash\gamma}\right)\ge c.
\end{equation}
\end{Proposition}

The proof of Proposition \ref{matrix:class} can be found in Section 5 (Appendix).

%format the remarks so that they do not indent.  I think you can do
%this with \noindent

% you need to put \noindent before each paragraph where you have a
% remark heading.

{\bf Remark 2.2.}\,\, Proposition \ref{matrix:class} demonstrates that 
Assumptions \ref{assump:2} and \ref{assump:8} can hold under general classes of design matrices.
One such class consists of matrices $X$ satisfying
\begin{equation}\label{matrix:requirement}
1/\bar{c}\le
\lambda_{-}\left(\frac{1}{n}X^TX\right)\le
\lambda_{+}\left(\frac{1}{n}X^TX\right)\le \bar{c},
\end{equation}
where $\bar{c}$ is some positive constant. For any $\gamma\in S_1$, 
we will have that $\gamma^0\subset\gamma$ and $\gamma^0\neq\gamma$. 
Thus, by Proposition \ref{matrix:class}, $\lambda_{-}\left(\frac{1}{n}X_{\gamma\backslash\gamma^0}'(I_n-P_{\gamma^0})X_{\gamma\backslash\gamma^0}\right)\ge 1/\bar{c}$,
i.e., inequality (\ref{cond_viii}) in Assumption \ref{assump:8} holds. 
Notice that when $\gamma\in S_2$, the relationship
$\gamma\subset\gamma^0\vee\gamma$ and $\gamma\neq\gamma^0\vee\gamma$ holds, 
where $\gamma^0\vee\gamma$ denotes the $p$-vector with $j$th component 
the larger of $(\gamma^0)_j$ and $\gamma_j$, then Assumption \ref{assump:2} 
follows by applying Proposition \ref{matrix:class}.\qed

In the following text, we assume that data are generated from the true
model $\textbf{y}=X\beta^0+\epsilon$ with $\epsilon\sim N(0,
\sigma_0^2 I_n)$. Let $\gamma^0$ be the $p$-dimensional state vector
corresponding to $\beta^0$. Unless otherwise stated, 
the limits in our main results will be taken when $n\rightarrow\infty$.

\begin{Theorem}\label{thm1} Suppose that $\gamma^0$ is nonnull and 
  Assumptions \ref{assump:1}--\ref{assump:4} and
  \ref{assump:6}--\ref{assump:8} are satisfied. Let $\delta\ge0$
  satisfy Assumption \ref{assump:8}. If
  $p^{\alpha_0}=o(n^{1-\delta}\underline{\phi}_n)$ for some $\alpha_0>2$, then
  $\sup\limits_{c_1,\ldots,c_p\in[\underline{\phi}_n,\bar{\phi}_n]}
  \max\limits_{\gamma\neq\gamma^0}\,\,p(\gamma|Z)/p(\gamma^0|Z)\rightarrow0$ in probability.
  If $p^{\alpha_0+2}=o(n^{1-\delta}\underline{\phi}_n)$ for some $\alpha_0>2$,
  then $\sup\limits_{c_1,\ldots,c_p\in[\underline{\phi}_n,\bar{\phi}_n]}
  \sum\limits_{\gamma\neq\gamma^0}p(\gamma|Z)\rightarrow0$ in probability,
  and consequently, $\inf\limits_{c_1,\ldots,c_p\in[\underline{\phi}_n,\bar{\phi}_n]}p(\gamma^0|Z)\rightarrow1$ in probability.
\end{Theorem}

The proof of Theorem \ref{thm1} follows
by first deriving asymptotic approximations of the posterior odds ratios
$p(\gamma|Z)/p(\gamma^0|Z)$ for any $\gamma\neq\gamma^0$, and then using
these approximations to show that $\sum\limits_{\gamma\neq\gamma^0}p(\gamma|Z)/p(\gamma^0|Z)\rightarrow 0$ in probability.
The limit $p(\gamma^0|Z)\rightarrow1$ (in probability) thus immediately follows from (\ref{PMC:BayesFactor}).
Details are in the Appendix.

{\bf Remark 2.3.}\,\, Theorem \ref{thm1} provides sufficient conditions under 
  which (\ref{mpp}) and PMC are satisfied. It asserts that, with large probability,
  uniformly for $c_j$'s $\in[\underline{\phi}_n,\bar{\phi}_n]$,
  $p(\gamma^0|Z)$ dominates $p(\gamma|Z)$ for any $\gamma\neq\gamma^0$, 
  and $p(\gamma^0|Z)$ approaches one in probability. Thus,
  with large probability, the true model $\gamma^0$ will be selected from a Bayesian perspective.\qed
  
{\bf Remark 2.4.}\,\,  When combined with certain dimension reduction techniques such as sure 
independence screening (SIS) proposed by Fan and Lv (2008),
one can generalize Theorem \ref{thm1} to the ultra-high dimensional setting, i.e., $p\gg n$.
This framework has been explored by many authors from non-Bayesian perspectives (see, e.g., 
Meinshausen and B\"{u}hlmann, 2006; Meinshausen and Yu, 2009; Zhang and Huang, 2010; B\"{u}hlmann and Kalisch, 2010).
Here, we explore it by a Bayesian way.
The basic idea is to first reduce the high-dimensional linear model so that the model dimension is below $n$,
and then apply Bayesian model (\ref{full-Bayesian-model}) to this reduced linear model.
Under suitable conditions and using the arguments similar to the proof of Theorem \ref{thm1}, 
one can show that the posterior probability of the true model based on the reduced linear model 
converges in probability to 1. We refer to Supplement A for the description of this 
result and details of the proof.\qed

The following result is an application of Theorem \ref{thm1} in a
special setting, which allows the growth rate of $p$ to be
$p\log{n}=o(n)$.

\begin{Corollary}\label{cor:3}
  Suppose that $\gamma^0$ is nonnull and Assumptions \ref{assump:1},
  \ref{assump:2} and inequality (\ref{cond_viii}) are
  satisfied. Assume that $\min\limits_{j\in\gamma^0}|\beta_j^0|\ge
  \psi_n$ with $\inf\limits_n\psi_n>0$, and $p$ satisfies
  $p\log{n}=o(n)$. Let 
$\delta\ge 0$ be as specified in
  inequality (\ref{cond_viii}) and
suppose
  there exists a constant $\delta_0$ with $\delta_0>3+\delta$ such
  that $k_n=O(n^{\delta_0})$.
Then with the selection $\bar{\phi}_n=O(n^{\delta_0})$
  and $n^{\delta_0}=O(\underline{\phi}_n)$, we have
  $\inf\limits_{c_1,\ldots,c_p\in[\underline{\phi}_n,\bar{\phi}_n]}p(\gamma^0|Z)\rightarrow1$ in probability.
\end{Corollary}

The proof of Corollary \ref{cor:3} can be finished by choosing
$\alpha_0 \in (2,\delta_0-\delta-1)$ and verifying the assumptions in
Theorem \ref{thm1}.

Theorem \ref{thm1} deals with the case when the true model is
nonnull. If the true model is null, then the response vector
$\textbf{y}$ will have a zero mean. The corresponding result is
summarized below.

\begin{Theorem}\label{thm3} Suppose $\gamma^0$ is null, i.e.,
  $\textbf{y}=\epsilon\sim N(0,\sigma_0^2 I_n)$, and that
  Assumptions \ref{assump:1} and
  \ref{assump:5}--\ref{assump:8} are satisfied. If
  $p^{\alpha_0}=o(n^{1-\delta}\underline{\phi}_n)$ for some $\alpha_0>2$, then
  $\sup\limits_{c_1,\ldots,c_p\in[\underline{\phi}_n,\bar{\phi}_n]}
  \max\limits_{\gamma\neq\gamma^0}\,\,p(\gamma|Z)/p(\gamma^0|Z)\rightarrow0$ in probability. If
  $p^{\alpha_0+2}=o(n^{1-\delta}\underline{\phi}_n)$ for some $\alpha_0>2$, then
  $\sup\limits_{c_1,\ldots,c_p\in[\underline{\phi}_n,\bar{\phi}_n]}
  \sum\limits_{\gamma\neq\gamma^0}p(\gamma|Z)\rightarrow 0$ in probability, and
  consequently, $\inf\limits_{c_1,\ldots,c_p\in[\underline{\phi}_n,\bar{\phi}_n]}p(\gamma^0|Z)\rightarrow1$ in probability.
\end{Theorem}

The proof of Theorem \ref{thm3} is similar to Theorem \ref{thm1} and can be found in
Appendix.

Although it is valid for a general type of design matrix, Theorem
\ref{thm1} requires that $p$ grows slower than $n$. More precisely, if
the Assumptions in Theorem \ref{thm1} are satisfied, then
$p=o(n)$. To see this, we notice that Assumptions
\ref{assump:6}, \ref{assump:7} and the fact that $\psi_n\le k_n^{1/2}$
lead to $\psi_n=O(n^{\delta_0})$ for some $\delta_0>0$. Therefore, $p=o(n)$
follows from Assumption \ref{assump:4}. In order to obtain consistency
when $p$ may grow as fast as $n$, one idea, but not the weakest
possible, is to assume orthogonality of $X$, i.e., $X^TX=nI_p$, and
to relax Assumption \ref{assump:7}. To simplify the technical proof,
we assume in the following Corollaries \ref{cor:1} and \ref{cor:2} that all 
$c_j$'s in model (\ref{full-Bayesian-model}) are equal to $\phi_n$.
Moreover, we need the following
assumption about the growth rates of $s_n$ and $p$ to replace
Assumptions \ref{assump:4} and \ref{assump:5}.

\begin{Assumption}\label{assump:9} Let $a_n=n+\sigma_0^{-2}
  k_n/(n^{-1}+\phi_n)$ and $\zeta\in (1,\infty)$ be a constant such that
  $n\psi_n^2>\sigma_0^2\zeta a_n$ as $n\rightarrow\infty$. 
  The numbers $p$ and $s_n$ with
  $p\rightarrow\infty$ and $s_n\le p\le n$ satisfy
\begin{enumerate}[(i).]
\item $s_n=o\left(\min\left\{\frac{(n+\nu)\log(n\psi_n^2/(\sigma_0^2\zeta a_n))}{\log(1+n\phi_n)},n\psi_n^2, n\right\}\right)$.
\item $p\log{p}=o\left(a_n\right)$. 
\end{enumerate}
\end{Assumption}

Assumption \ref{assump:9} potentially allows the case $p=n$. To see
this, suppose $s_n=O(1)$ and we choose $\phi_n$ such that $(n+\nu)/\log(1+n\phi_n)\rightarrow\infty$.
When $a_n$ grows faster than $n\log{n}$ and 
$n\psi_n^2/a_n\rightarrow\infty$, $p=n$ will satisfy Assumption
\ref{assump:9}. However, this requires $\psi_n^2$ to grow at least faster
than $\log{n}$. This extra requirement on $\psi_n^2$ has not been
imposed by Theorems \ref{thm1} and \ref{thm3}, and can be treated as
the price which we pay to relax the growth rate for $p$. Under
Assumption \ref{assump:9} and assuming orthogonality on $X$, we have
the following consistency result which allows a faster
growth rate for the dimension $p$.

\begin{Corollary}\label{cor:1}
  Assume that $X^TX=nI_p$ and $\Sigma=\phi_n I_p$ with
  $n\phi_n\rightarrow\infty$ and $\log{\phi_n}=O(\log{n})$. 
  Suppose $\gamma^0$ is nonnull and that
  Assumptions \ref{assump:1} and
  \ref{assump:9} are satisfied. If
  $p^{\alpha_0(n+\nu)/a_n}=o(n\phi_n)$ for some $\alpha_0>2$, then
  $\max\limits_{\gamma\neq\gamma^0}\,\,p(\gamma|Z)/p(\gamma^0|Z)\rightarrow0$ in probability.
  If $p=o\left((n+\nu)\log\left(\frac{n\psi_n^2}{\sigma_0^2\zeta a_n}\right)\right)$
  with $\zeta$ specified in Assumption \ref{assump:9}, 
  and $p^{2+\alpha_0(n+\nu)/a_n}=o(n\phi_n)$ for some $\alpha_0>2$,
  then $\sum\limits_{\gamma\neq\gamma^0}p(\gamma|Z)\rightarrow0$ in probability,
  and consequently, $p(\gamma^0|Z)\rightarrow1$ in probability.
\end{Corollary}

The proof of Corollary \ref{cor:1} is similar to those for Theorems
\ref{thm1} and \ref{thm3} and is given in Supplement B. The
following result, which requires a special model set-up, 
demonstrates that PMC and consistency of the posterior odds ratio may
hold in some situations but fail in others.

\begin{Corollary}\label{cor:2}
  Assume $p=n$, $X^TX=nI_n$ and $\Sigma=\phi_n I_n$. Suppose
  $\min\limits_{j\in\gamma^0}|\beta_j^0|\ge\psi_n$ with
  $\psi_n^2=c_1n^{1+\delta_1}(\log{n})^2$ for some constants
  $\delta_1>1$ and $c_1>0$, $k_n=O(\psi_n^2)$ and
  $p(\gamma)=\textrm{constant}$ for all $\gamma$. Assume that $s_n=s$
  with $s>0$ a fixed integer, i.e., the true parameter vector $\beta^0$
  contains exactly $s$ nonzero components.
\begin{enumerate}[(a).]
\item Suppose $\phi_n=c_2n^{\delta_2}$ for some constants $c_2>0$ and $\delta_2$. 
\begin{enumerate}[i.]
\item If $-1<\delta_2\le 1$, then
  $\max\limits_{\gamma\neq\gamma^0}\,\,p(\gamma|Z)/p(\gamma^0|Z)\rightarrow0$ in probability,
  but PMC does not hold. Specifically, when $-1<\delta_2<1$,
  $p(\gamma^0|Z)\rightarrow0$, a.s.; when $\delta_2=1$, then there
  exists a constant $c_0$ with $0<c_0<1$ such that $\limsup\limits_n
  p(\gamma^0|Z)\le c_0$, a.s.
\item If $1<\delta_2\le \delta_1$, then $p(\gamma^0|Z)\rightarrow1$ in probability.
\end{enumerate}
\item If $n^{n\log{n}}=O(\phi_n)$, then
  $p(\emptyset|Z)/p(\gamma^0|Z)\rightarrow\infty$ in probability, where $\emptyset$
  represents the null model. Therefore, $p(\gamma^0|Z)\rightarrow0$ in probability.
\item If $n\phi_n\rightarrow \eta\in [0,\infty)$, then almost surely,
  $\liminf\limits_n
  \max\limits_{\gamma\neq\gamma^0}\,\,p(\gamma|Z)/p(\gamma^0|Z)\ge
  (1+\eta)^{-1/2}$ and $\lim\limits_{n} p(\gamma^0|Z)=0$.
\end{enumerate}
\end{Corollary}

The proof of Corollary \ref{cor:2} is given in Supplement B.

{\bf Remark 2.5.}\,\, The main contribution of Corollary \ref{cor:2}
is to demonstrate the difference between PMC and (\ref{mpp}), and
provide example growth rates for $\phi_n$ under which the two forms of
consistency fail. Although this is obtained in a special situation,
similar results should be still true under a more general setting, for
instance, where $p<n$ or $X^TX$ is not diagonal, but we do not
  consider those circumstances here.

Corollary \ref{cor:2} (a) demonstrates that (\ref{mpp}) does not
necessarily imply PMC. This means that, although the posterior
probability of the true model might not be approaching one, the ratio
of the posterior probabilities of any ``incorrect" model and the true
model can still converge to zero. This phenomenon will not occur when
$p$ is fixed. In practice, (\ref{mpp}) is sufficient to make a correct
model selection even if PMC might fail.

Corollary \ref{cor:2} (b) and (c) demonstrate that in order to make a
correct model selection, $\phi_n$ cannot be either too small or too large. 
Specifically, when $\phi_n=o(n^{-1})$, it follows by Corollary \ref{cor:2} (c) that
almost surely $\liminf\limits_n
\max\limits_{\gamma\neq\gamma^0}\,\,p(\gamma|Z)/p(\gamma^0|Z)\ge
1$. Thus, with probability one, for any $\varepsilon>0$,
there exists an integer $N$ such that for any $n\ge N$
\[
\max\limits_{\gamma\neq\gamma^0}\,\,p(\gamma|Z)/p(\gamma^0|Z)\ge 1-\varepsilon.
\]
This implies that there exists a model, say $\gamma^*$, such that $p(\gamma^*|Z)\ge (1-\varepsilon)p(\gamma^0|Z)$.
Thus, when $\varepsilon$ is small, either $p(\gamma^*|Z)>p(\gamma^0|Z)$, 
or $p(\gamma^*|Z)$ is very close to $p(\gamma^0|Z)$,
which will both affect the selection result. On the other hand, when $\phi_n$ is growing
faster than $n^{n\log{n}}$, it follows from (b) that the null model will
be preferred in favor of $\gamma^0$.

Corollary \ref{cor:2} (b) and (c) can be also understood intuitively.
When $\phi_n$ is too small, the two distribution components in the mixture prior of $\beta$
tend to be indistinguishable so that it is difficult to separate the true model from some
incorrect model; when $\phi_n$ approaches infinity, by (\ref{eq:post}), 
the posterior probability of any nonnull model approaches zero, and thus, all 
$\beta_j$'s are forced to be zero. This conclusion has been empirically
obtained by Smith and Kohn (1996) under spline regression models.\qed

{\bf Remark 2.6.}\,\, Using arguments similar to the proofs of Theorems \ref{thm1} and \ref{thm3},
and by the Borel-Cantelli lemma (see Shao, 2003),
one can show the almost sure convergence of $p(\gamma^0|Z)$. We refer to Supplement C
for details.\qed

To conclude this section, let us
look at an example which demonstrates that, when $\bar{\phi}_n=\bar{\phi}$ and $\underline{\phi}_n=\underline{\phi}$ 
with $\bar{\phi}$ and $\underline{\phi}$ unrelated to $n$, 
consistency might still hold under certain circumstances. 
This is motivated by a full Bayesian framework which requires all hyperparameters 
to be fixed.

{\bf Example 2.1.}\,\, If a full Bayesian approach is desired, then we
have to preselect the hyperparameters $c_j$'s, and so $\bar{\phi}_n=\bar{\phi}$ and
$\underline{\phi}_n=\underline{\phi}$ could be
fixed. Assume that $k_n=O(1)$, which is a slightly weaker
assumption than that in Jiang (2007). Note that Assumptions
\ref{assump:6} and \ref{assump:7} follow immediately. Suppose
$\min\limits_{j\in\gamma^0}|\beta_j^0|\ge\psi_n$ with $\psi_n\propto
n^{-1/4}\sqrt{\log{n}}$, the prior distribution of model $\gamma$
satisfies Assumption \ref{assump:1}. Assume that $s_n=s$ with $s>0$ 
a fixed integer (thus, the true model is nonnull), and design
matrix $X$ satisfies (\ref{matrix:requirement}). Therefore, by Proposition \ref{matrix:class}
and Remark 2.2, Assumptions \ref{assump:2} and 
\ref{assump:8} both hold. We also notice that Assumption \ref{assump:3} is well
satisfied. It follows from Theorem \ref{thm1} that if $p\propto n^r$ for
some $0<r<1/2$, then with probability approaching one, (\ref{mpp}) holds,
i.e., the true model can be correctly selected; if $p\propto n^r$
for some $0<r<1/4$, then PMC holds in probability.

\section{Generalizations to $g$-prior settings}

In section 2, we assume in the Bayesian model
(\ref{full-Bayesian-model}) that the prior variance of a nonzero
$\beta_j$ is $c_j\sigma^2$ with $c_j$ being  fixed a priori. In
practice, one may consider placing a prior distribution $g(c)$ on the
$c_j$'s, which reduces to the so-called $g$-prior setting (see Zellner
1986; Liang \textit{et al.} 2008). In this section, we will give some
asymptotic results under a $g$-prior setting.

We consider the following variation in model
(\ref{full-Bayesian-model}):
\begin{eqnarray*}
  \beta_j|\gamma_j,\sigma^2, c&\sim& (1-\gamma_j)\delta_0+\gamma_j N(0,c\sigma^2),\,\, j=1,\ldots,p,\nonumber\\
  c&\sim& g(c),
\end{eqnarray*}
where $g$ is a proper prior distribution on $[0,\infty)$. We still use
$p(\gamma^0|Z)$ to denote the posterior probability of the true
model. Note that $p(\gamma|Z)$ is obtained by integrating
$p(\gamma,\beta,\sigma^2,c|Z)$ with respect to $(\beta,\sigma^2,c)$.

\begin{Theorem}\label{thm:g-prior:1}
  Suppose that $\gamma^0$ is nonnull and Assumption \ref{assump:1}
  holds.  Furthermore, suppose $\|\beta^0\|_2=O(1)$, $s_n=O(1)$,
  $\min\limits_{j\in \gamma^0}|\beta^0_j|\ge \psi_n$ with
  $\psi_n\propto n^{-1/4} \sqrt{\log{n}}$, and the design matrix $X$
  satisfies property (\ref{matrix:requirement}).
\begin{enumerate}[(i)]
\item Let the support of $g$ be $[\underline{\phi},\bar{\phi}]$ with
  $0<\underline{\phi}<\bar{\phi}<\infty$.  If $p\propto n^{r}$ for
  some $0<r<1/2$, then $\max\limits_{\gamma\neq\gamma^0}
  p(\gamma|Z)/p(\gamma^0|Z)\rightarrow 0$ in probability.
\item Let $g$ be proper on $[0,\infty)$. If $p\propto n^{r}$ for some
  $0<r<1/4$, then $p(\gamma^0|Z)\rightarrow 1$ in probability.
\end{enumerate} 
\end{Theorem}  

\begin{Theorem}\label{thm:g-prior:2}
  Suppose that $\gamma^0$ is nonnull and Assumptions
    \ref{assump:1}, \ref{assump:3} and \ref{assump:4} are satisfied.
    Let $X$ satisfy (\ref{matrix:requirement}). Suppose that
    $\underline{\phi}_n$ and $\bar{\phi}_n$ satisfy Assumptions
    \ref{assump:6} and \ref{assump:7}.
\begin{enumerate}[(i)]
\item Let the support of $g$ be $[\underline{\phi}_n,\bar{\phi}_n]$. If
  $p^{\alpha_0}=o(n\underline{\phi}_n)$ for some $\alpha_0>2$, then
  $\max\limits_{\gamma\neq\gamma^0}
  p(\gamma|Z)/p(\gamma^0|Z)\rightarrow 0$ in probability.
\item Let $g$ be proper on $[0,\infty)$ such that
  $1-\int_{\underline{\phi}_n}^{\bar{\phi}_n} g(c)dc=o(1)$.  If
  $p^{\alpha_0+2}=o(n\underline{\phi}_n)$ for some $\alpha_0>2$, then
  $p(\gamma^0|Z)\rightarrow 1$ in probability.
\end{enumerate}
\end{Theorem}

The proof of Theorems \ref{thm:g-prior:1} and
  \ref{thm:g-prior:2} is given in the Appendix.
  
{\bf Remark 3.1.}\,\, Theorems \ref{thm:g-prior:1} and
  \ref{thm:g-prior:2} provide sufficient conditions for (\ref{mpp})
  and PMC under a $g$-prior setting. They state that with large
  probability, $p(\gamma^0|Z)$ dominates $p(\gamma|Z)$ for any
  $\gamma\neq\gamma^0$, and $p(\gamma^0|Z)$ approaches one in
  probability.  In particular, the prior $g$ in Theorem
  \ref{thm:g-prior:1} does not depend on $n$, which corresponds to a
  full Bayesian framework, but we need to impose a narrow restriction on the
  growth rate of $p$, namely, that $p$ is growing slower than $n^{1/4}$ or $n^{1/2}$,
  corresponding to PMC or (\ref{mpp}). In Theorem \ref{thm:g-prior:2}  
  $g$ might depend
  on $n$, but we can allow $p$ to grow faster with $n$.  \qed

{\bf Remark 3.2.}\,\, We conjecture, although do not rigorously prove, that the ranges
$0<r<1/2$ and $0<r<1/4$ in parts (a) and (b) of Theorems \ref{thm:g-prior:1} are optimal,
in the sense that for any $r>1/2$, if $p\propto n^r$, then $\max\limits_{\gamma\neq\gamma^0} p(\gamma|Z)/p(\gamma^0|Z)$ 
does not converge to zero in probability; and for any $r>1/4$, if $p\propto n^r$,
then $p(\gamma^0|Z)$ does not converge to one in probability.

{\bf Remark 3.3.}\,\, Liang {\it et al.} (2008) obtained model consistency under a mixture $g$-prior setting.
Their proof relies on the Laplace approximation of the integrals.
While the proofs of both Theorems \ref{thm:g-prior:1} and \ref{thm:g-prior:2} rely
on the uniform convergence in Theorem \ref{thm1}.\qed  
\section{Numerical results}

This paper has been concerned with asymptotic properties of Bayesian
posterior probabilities. In this section, we briefly explore the 
finite
sample behavior of the model selection procedure for a few different
prior settings and different rates of growth for $p$. Our basic
approach is to simulate observations from model (\ref{basic:linear:model}), employ the model
selection process, and summarize the results.

To construct random design matrices $X$, we generated \textit{iid}
$p$-dimensional row vectors $U_1,\ldots, U_n\sim N(\textbf{0}, I_p)$
and let $U$ be an $n\times p$ matrix with $i$th row $U_i$ for
$i=1,\ldots,n$. Then we let $X=\sqrt{n}U \left(U^TU\right)^{-1/2}$.
Thus, $X^TX=nI_p$. (We choose $X$ to be orthonormal for purposes of
illustration, although, as we saw in the preceding material, results
can be derived for general $X$.)  To explore the dimension effect, we
have considered three growth rates for $p$ with respect to $n$: (1)
$p=n^{1/4}$, (2) $p=n^{1/2}$ and (3) $p=n^{3/4}$. Data were simulated
from model (\ref{model}) with $\sigma=1$, $s_n=2$ and the true model
coefficients $(\beta_1^0,\beta_{2}^0)=(2,2)$ and
$(\beta_{3}^0,\ldots,\beta_{p}^0)=(0,\ldots,0)$.  We considered sample
sizes $n=100$, $200$ and $400$ respectively.

The hierarchical Bayesian model (\ref{full-Bayesian-model}) was fitted
and the prior distributions on $\sigma^2$ and $\gamma$ were assumed to
be $1/\sigma^2\sim\chi_4^2$ and $p(\gamma_j=1)=w_j$, for any
$j=1,\ldots,p$. We examined two cases for the $w_j$'s, namely, Case I:
$w_j=0.5$ for $1\le j\le p$; and Case II: $w_1=w_2=0.3$,
$w_3=\ldots=w_p=0.7$.  Case I places equal prior probabilities on all
the models, while Case II places larger prior probabilities on the
``incorrect" models.  For simplicity, we let
$c_1=\ldots=c_p=\phi_n$. The values of $\phi_n$ were chosen to be
$\phi_n=10, 100, 1000$. A total of 20,000 samples of
$(\beta,\gamma,\sigma)$ were drawn from the posterior distribution
$p(\beta,\gamma,\sigma|Z)$ using a sub-blockwise Gibbs sampler
developed by Godsill and Rayner (1998). We recorded the last $10,000$
samples and treated the previous 10,000 samples as burnins.
Convergence was assessed by applying Gelman-Rubin's statistic to
5 parallel Markov chains for each $\phi_n$. If we denote
$\gamma^{(1)},\ldots,\gamma^{(10000)}$ to be the last $10,000$ samples
of $\gamma$, then $p(\gamma^0|Z)$ is approximated by
$p(\gamma^0|Z)\approx \sum\limits_{t=1}^{10000}
I(\gamma^{(t)}=\gamma^0)/10000.$

To study the frequentist behavior of $p(\gamma^0|Z)$, we have
generated 100 data sets $Z_1,\ldots,Z_{100}$ independently from model
(\ref{model}), and for each $\phi_n$ calculated the corresponding 100
posterior probabilities $p(\gamma^0|Z_m)$, $m=1,\ldots,100$ as
described in the preceding paragraph. This idea was inspired from
Fern\'{a}ndez {\it et al.} (2001) who studied the Bayesian selection
problem when $p$ is fixed.

Table \ref{Table:1} summaries the mean and standard deviations of the
100 $p(\gamma^0|Z_m)$'s.  We compared four settings. Specifically, Setting 1
to 3 correspond to $\phi_n=10, 100, 1000$ under the Bayesian model
(\ref{full-Bayesian-model}), and Setting 4 uses a hyper $g$-prior with
tuning parameter 3 (see Liang {\it et al.} 2008).  Setting 4 was
performed using the R package BAS available from
\textit{http://www.stat.duke.edu/$\sim$clyde/BAS}.  We observe that
when $p=n^{1/4}$, all four settings select the true model with high
posterior probability. For the faster growth rates $p=n^{1/2}$ and
$p=n^{3/4}$, the results are more mixed. Generally,  Setting
1 performs the worst and Setting 3 performs the best. In summary,
when $p$ is small compared to $n$, fixing $\phi_n$ to be 10, 100 or 1000
will result in equally good results; when $p$ is larger (compared to $n$),
we recommend using $\phi_n=1000$ under model
(\ref{full-Bayesian-model}), if good asymptotics behavior is of
interest.
%are these results consistent with the asymptotics?
%comment on the change over n

\begin{table}[htp]
\begin{center}
\begin{tabular}{ccccccccc}
                            &                         &&\multicolumn{2}{c}{$n=100$}&\multicolumn{2}{c}{$n=200$}&\multicolumn{2}{c}{$n=400$}\\ 
                            &                         &          & mean & std  & mean & std  & mean & std \\ \hline
\multirow{8}{*}{$p=n^{1/4}$}&\multirow{4}{*}{Case I}  & Setting 1 & 0.94 & 0.05 & 0.96 & 0.04 & 0.92 & 0.10\\
                            &                         & Setting 2 & 0.98 & 0.02 & 0.99 & 0.02 & 0.97 & 0.05\\
                            &                         & Setting 3 & 0.99 & 0.01 & 0.99 & 0.01 & 0.99 & 0.02\\
                            &                         & Setting 4 & 0.96 & 0.04 & 0.96 & 0.05 & 0.95 & 0.05\\ \cline{2-9}
                            & \multirow{4}{*}{Case II}& Setting 1 & 0.86 & 0.14 & 0.91 & 0.08 & 0.86 & 0.10\\
                            &                         & Setting 2 & 0.94 & 0.10 & 0.97 & 0.04 & 0.95 & 0.05\\
                            &                         & Setting 3 & 0.98 & 0.05 & 0.99 & 0.02 & 0.98 & 0.02\\
                            &                         & Setting 4 & 0.92 & 0.05 & 0.94 & 0.05 & 0.89 & 0.08\\ \hline 
\multirow{8}{*}{$p=n^{1/2}$}&\multirow{4}{*}{Case I}  & Setting 1 & 0.60 & 0.14 & 0.56 & 0.14 & 0.53 & 0.12\\
                            &                         & Setting 2 & 0.82 & 0.10 & 0.81 & 0.11 & 0.80 & 0.10\\
                            &                         & Setting 3 & 0.94 & 0.05 & 0.93 & 0.06 & 0.93 & 0.05\\
                            &                         & Setting 4 & 0.68 & 0.10 & 0.63 & 0.13 & 0.62 & 0.12\\ \cline{2-9}
                            &\multirow{4}{*}{Case II} & Setting 1 & 0.34 & 0.12 & 0.29 & 0.11 & 0.27 & 0.11\\
                            &                         & Setting 2 & 0.65 & 0.14 & 0.63 & 0.14 & 0.62 & 0.14\\
                            &                         & Setting 3 & 0.86 & 0.09 & 0.85 & 0.10 & 0.84 & 0.09\\
                            &                         & Setting 4 & 0.42 & 0.10 & 0.41 & 0.10 & 0.36 & 0.10\\ \hline
\multirow{8}{*}{$p=n^{3/4}$}&\multirow{4}{*}{Case I}  & Setting 1 & 0.14 & 0.07 & 0.07 & 0.04 & 0.04 & 0.03\\
                            &                         & Setting 2 & 0.47 & 0.13 & 0.38 & 0.12 & 0.33 & 0.10\\
                            &                         & Setting 3 & 0.77 & 0.10 & 0.71 & 0.13 & 0.68 & 0.11\\
                            &                         & Setting 4 & 0.21 & 0.08 & 0.16 & 0.05 & 0.16 & 0.06\\ \cline{2-9}
                            &\multirow{4}{*}{Case II} & Setting 1 & 0.02 & 0.01 & 0.00 & 0.00 & 0.00 & 0.00\\
                            &                         & Setting 2 & 0.20 & 0.10 & 0.13 & 0.06 & 0.08 & 0.05\\
                            &                         & Setting 3 & 0.55 & 0.15 & 0.48 & 0.12 & 0.41 & 0.12\\
                            &                         & Setting 4 & 0.04 & 0.02 & 0.03 & 0.02 & 0.04 & 0.02\\ \hline                                        
\end{tabular}
\end{center}
\caption{\textit{Means and standard deviations of the 100 $p(\gamma^0|Z_m)$'s.
Settings 1 to 3 correspond to $\phi_n=10, 100, 1000$ under the Bayesian model 
(\ref{full-Bayesian-model}), and Setting 4 uses hyper $g$-prior with tuning parameter 3.}}
\label{Table:1}
\end{table}

\section{Conclusion}

Previous work about posterior model consistency (PMC) includes
Fern\'{a}ndez {\it et al.} (2001) and Liang {\it et al.} (2008) when
the number of parameters $p$ is fixed. In this paper, we have studied
PMC when the model dimension $p$ grows with sample size
$n$. Specifically, we have shown that, under a variation of the
Bayesian model proposed by George and McCulloch (1993), the posterior
probability of the true model converges to one, i.e., PMC holds. We
have obtained this result in two situations: (i) the design matrix $X$
is general while $p$ grows slower than $n$, e.g., $p\log{n}=o(n)$;
(ii) $X^TX/n$ is the identity matrix and $p$ may grow as fast as $n$,
e.g., $p=n$. Furthermore, we have demonstrated under a special
framework that the consistency results may fail if $\phi_n$ is too
small or too large, where $\phi_n$ is the hyperparameter controlling
the prior variance of the nonzero model coefficients.  More precisely,
when $\phi_n=o(n^{-1})$ (an example of small order) or when
$n^{n\log{n}}=O(\phi_n)$ (an example of large order), both PMC and
consistency of the posterior odds ratio fail. Besides that, our
results do not require that the candidate models are pairwise nested.

Berger {\it et al.} (2003), Moreno {\it et al.} (2010) 
and Gir\'{o}n {\it et al.} (2010) have proved the consistency of Bayes factor
when $p$ is growing with $n$. This form of
consistency, under our framework, is equivalent to the consistency of
the posterior odds ratio if the prior odds ratio is uniformly bounded
from above and below, so it is of interest to illustrate the
relationship between PMC and consistency of posterior odds ratio. We
have considered a special framework and shown that PMC implies
consistency of the posterior odds ratio but the reverse may not be
true. This is different from the finding by Liang {\it et al.} (2008) who
demonstrate the equivalence of PMC and consistency of the Bayes factor
when $p$ is fixed. When combined with dimension reduction procedures such as SIS (Fan and Lv, 2008), 
our results can be also extended to ultrahigh-dimensional situations. 
We have also generalized the consistency results to a $g$-prior setting studied by Zellner (1986)
and Liang {\it et al.} (2008).

We close with an observation about extending the current results.
Assumption \ref{assump:7} is a technical assumption
used to facilitate the proof and may not be
the weakest possible. We leave it to future work to determine whether this condition
can be further weakened or even removed. 

\section{Appendix: proofs} 

In this section, we prove the main results in Section 2.
We also prove some lemmas which are useful to establish the main results.
Let $pr(\cdot)$ denote the probability measure associated with the underlying
probability space.

{\bf Proof of Proposition \ref{matrix:class}.}\,\, It follows by 
assumption that $\frac{1}{n}X^T_{\bar{\gamma}}X_{\bar{\gamma}}\ge cI_{|\bar{\gamma}|}$. 
Letting $X_{\bar{\gamma}}=\left(X_\gamma,X_{\bar{\gamma}\backslash\gamma}\right)$, we can write 
$\frac{1}{n}X^T_{\bar{\gamma}}X_{\bar{\gamma}}=
\left(\begin{array}{cc}A&B \\ B^T&C\end{array}\right)$,
where $A=X_{\gamma}^TX_{\gamma}/n$, $B=X_{\gamma}^T X_{\bar{\gamma}\backslash\gamma}/n$ and $C=X_{\bar{\gamma}\backslash\gamma}^TX_{\bar{\gamma}\backslash\gamma}/n$. 
By formula for the inverse of blocked matrix (Seber and Lee, 2003, page 466),
the lower right corner of $\left(\frac{1}{n}X^T_{\bar{\gamma}}X_{\bar{\gamma}}\right)^{-1}$
is $B_{22}^{-1}$ with 
$B_{22}=C-B^TA^{-1}B=\frac{1}{n}X_{\bar{\gamma}\backslash\gamma}^T(I_n-P_{\gamma})X_{\bar{\gamma}\backslash\gamma}$.
Then $B_{22}^{-1}\le c^{-1} I$, which implies $\lambda_{-}(B_{22})\ge c$. \qed

\begin{Lemma}\label{lemma1} Suppose $\epsilon\sim N(0,\sigma_0^2 I_n)$. Then:
\begin{enumerate}[(a).]

\item Let $v_\gamma=(I_n-P_\gamma)X_{\gamma^0\backslash\gamma}\beta^0_{\gamma^0\backslash\gamma}$. If $S_2$ is nonnull, then $\max\limits_{\gamma\in S_2} |v_\gamma^T \epsilon|/ \|v_\gamma\|_2=O_p\left(\sqrt{p}\right)$, where we adopt the convention that $|v_\gamma^T \epsilon|/\|v_\gamma\|_2=0$ when $v_\gamma=0$. 

\item If $S_1$ is nonnull, then for any $\alpha>2$, with probability approaching one,
$\max\limits_{\gamma\in S_1}\epsilon^T (P_\gamma-P_{\gamma^0})\epsilon/(|\gamma|-s_n) \le \alpha\sigma_0^2\log{p}$.

\item If $S_2$ is nonnull, and we adopt the convention that $\epsilon^T P_\gamma\epsilon/|\gamma|=0$ when $\gamma$ is
null, then for any $\alpha>2$, with probability approaching one,
$\max\limits_{\gamma\in S_2}\epsilon^T P_\gamma\epsilon/|\gamma| \le \alpha\sigma_0^2\log{p}$.

\end{enumerate}
\end{Lemma}

{\bf Proof of Lemma \ref{lemma1}.}\,\, We prove the result for the case where $X$ is deterministic,
and briefly talk about the proofs for the case where $X$ is random and independent of $\epsilon$.

(a) We first assume that $X$ is deterministic. By inequality (9.3) in Durrett (2005), 
if $\xi\sim N(0,1)$, then there exists a $C_0$ such that for any $t>1$, 
$pr(|\xi|\ge t)\le C_0 \exp\left(-t^2/2\right)$. 
Note that $|v_\gamma^T \epsilon|/(\sigma_0\|v_\gamma\|_2)\sim N(0,1)$, 
and therefore, by Bonferroni's inequality,
\[
pr\left(\max\limits_{\gamma\in S_2}\frac{|v_\gamma^T \epsilon|}{\|v_\gamma\|_2}\ge t\right)\le\sum\limits_{\gamma\in S_2} pr\left(\frac{|v_\gamma^T \epsilon|}{\|v_\gamma\|_2}\ge t\right)\le C_0 2^p \exp\left(-\frac{t^2}{2\sigma_0^2}\right).
\]
Then the result holds by setting $t=C\sigma_0\sqrt{2p}$ with large $C$. When $X$ is random but independent of $\epsilon$,
note that the conditional distribution of $|v_\gamma^T \epsilon|/(\sigma_0\|v_\gamma\|_2)$ given $X$ is $N(0,1)$. 
Thus, the proof can be finished by the above arguments.

(b) Suppose $X$ is deterministic. First, if $\xi=\chi_\mu^2$, then by Chebyshev's inequality, for any
$2<\alpha'<\alpha$, 
\begin{eqnarray*}
&&pr(\xi\ge \alpha\mu\log{p})\\&=&pr\left(\exp(\xi/\alpha')\ge \exp((\alpha/\alpha')\mu\log{p})\right)\\
&\le& \exp(-(\alpha/\alpha')\mu\log{p}) E\left\{\exp(\xi/\alpha')\right\}\\
&=& (1-2/\alpha')^{-\mu/2} \exp(-(\alpha/\alpha')\mu\log{p}).
\end{eqnarray*} Using this inequality, Bonferroni's inequality, and the fact that when $\gamma\in S_1$, $\epsilon^T(P_\gamma-P_{\gamma^0})\epsilon\sim\sigma_0^2\chi_{|\gamma|-s_n}^2$, we have
\begin{eqnarray*}
&&pr\left(\max\limits_{\gamma\in S_1}\frac{\epsilon^T (P_\gamma-P_{\gamma^0})\epsilon}{|\gamma|-s_n} \ge \alpha\sigma_0^2\log{p}\right)\\
&\le&\sum\limits_{\gamma\in S_1}pr\left(\epsilon^T (P_\gamma-P_{\gamma^0})\epsilon\ge \alpha\sigma_0^2 (|\gamma|-s_n)\log{p}\right)\\
&\le&\sum\limits_{\gamma\in S_1}(1-2/\alpha')^{-(|\gamma|-s_n)/2} \exp(-(\alpha/\alpha')(|\gamma|-s_n)\log{p})\\
&=&\sum\limits_{r=1}^{p-s_n}{p-s_n\choose r} (1-2/\alpha')^{-r/2}
\exp(-(\alpha/\alpha')r\log{p})\\
&=&\left(1+(1-2/\alpha')^{-1/2}p^{-\alpha/\alpha'}\right)^{p-s_n}-1\rightarrow0.
\end{eqnarray*}
When $X$ is random and independent of $\epsilon$, then conditioning on $X$, $\epsilon^T(P_\gamma-P_{\gamma^0})\epsilon\sim\sigma_0^2\chi_{|\gamma|-s_n}^2$.
Thus, the conclusion follows from the above arguments.

(c) We let $X$ be deterministic. The case where $X$ is random can be handled similarly.
Assume that $S_2$ contains nonnull models, and note that when $\gamma$ is nonnull, 
$\epsilon^TP_\gamma\epsilon\sim\sigma_0^2\chi_{|\gamma|}^2$. 
Fix arbitrarily $\alpha'$ such that $2<\alpha'<\alpha$. 
Then by the proof of part (b) we have
\begin{eqnarray*}
&&pr\left(\max\limits_{\gamma\in S_2}\frac{\epsilon^T P_\gamma\epsilon}{|\gamma|}\ge \alpha\sigma_0^2\log{p}\right)\\
&=&pr\left(\max\limits_{\gamma\in
S_2\backslash\{\emptyset\}}\frac{\epsilon^T P_\gamma\epsilon}{|\gamma|} \ge \alpha\sigma_0^2\log{p}\right)\\
&\le&\sum\limits_{\gamma\in S_2\backslash\{\emptyset\}}pr\left(\epsilon^T P_\gamma\epsilon\ge \alpha\sigma_0^2|\gamma|\log{p}\right)\\
&\le&\sum\limits_{\gamma\in S_2\backslash\{\emptyset\}} (1-\alpha'/2)^{-|\gamma|/2}\exp(-(\alpha/\alpha')|\gamma|\log{p})\\
&\le&\sum\limits_{r=1}^{p} {p\choose r} (1-2/\alpha')^{-r/2}p^{-(\alpha/\alpha')r}\\
&=&\left(1+(1-2/\alpha')^{-1/2}p^{-\alpha/\alpha'}\right)^{p}-1\rightarrow0.\qed
\end{eqnarray*}

{\bf Proof of Theorem \ref{thm1}.}\,\, We have
\begin{eqnarray}\label{eq5}
-\log{\left(p(\gamma|Z)/p(\gamma^0|Z)\right)}&=&-\log{\left(\frac{p(\gamma)}{p(\gamma^0)}\right)}+\frac{1}{2}\log{\left(\frac{\det(W_\gamma)}{\det(W_{\gamma^0})}\right)}\nonumber\\
&&+\frac{n+\nu}{2}\log{\left(\frac{1+\textbf{y}^T(I_n-X_\gamma U_\gamma^{-1}
X_\gamma^T)\textbf{y}}{1+\textbf{y}^T(I_n-X_{\gamma^0} U_{\gamma^0}^{-1} X_{\gamma^0}^T )\textbf{y}}\right)}\nonumber\\
&=&-\log{\left(\frac{p(\gamma)}{p(\gamma^0)}\right)}+ \frac{1}{2}\log{\left(\frac{\det(W_\gamma)}{\det(W_{\gamma^0})}\right)}\nonumber\\
&&+\frac{n+\nu}{2}\log{\left(\frac{1+\textbf{y}^T(I_n-X_\gamma U_\gamma^{-1}
X_\gamma^T)\textbf{y}}{1+\textbf{y}^T(I_n-P_\gamma)\textbf{y}}\right)}\nonumber\\
&&-\frac{n+\nu}{2}\log{\left(\frac{1+\textbf{y}^T(I_n-X_{\gamma^0} U_{\gamma^0}^{-1} X_{\gamma^0}^T)\textbf{y}}{1+\textbf{y}^T(I_n-P_{\gamma^0})\textbf{y}}\right)}\nonumber\\
&&+\frac{n+\nu}{2}\log{\left(\frac{1+\textbf{y}^T(I_n-P_\gamma)\textbf{y}}{1+\textbf{y}^T(I_n-P_{\gamma^0})\textbf{y}}\right)}.
\end{eqnarray}
Denote the above summands by $T_1, T_2, T_3, T_4, T_5$. 
By Assumption  \ref{assump:6}, $T_1$ is bounded below. 
Since $U_\gamma \ge X_\gamma^T X_\gamma$, we have $T_3\ge0$ for any $n$.

To approximate $T_4$, let 
\[
\Delta=\textbf{y}^T X_{\gamma^0}\left(X_{\gamma^0}^T X_{\gamma^0}\right)^{-1} \left(\Sigma_{\gamma^0}+\left(X_{\gamma^0}^T X_{\gamma^0}\right)^{-1}\right)^{-1} \left(X_{\gamma^0}^T X_{\gamma^0}\right)^{-1} X_{\gamma^0}^T \textbf{y}. 
\]
By the Sherman-Morrison-Woodbury matrix identity (Seber and Lee, 2003,page 467),
\begin{equation}\label{SMW:formula}
U_{\gamma^0}^{-1}-\left(X_{\gamma^0}^T X_{\gamma^0}\right)^{-1}=-\left(X_{\gamma^0}^T X_{\gamma^0}\right)^{-1}\left(\Sigma_{\gamma^0}+\left(X_{\gamma^0}^T X_{\gamma^0}\right)^{-1}\right)^{-1}\left(X_{\gamma^0}^T X_{\gamma^0}\right)^{-1}.
\end{equation}
By (\ref{SMW:formula}) and the fact that $\left(\Sigma_{\gamma^0}+\left(X_{\gamma^0}^T X_{\gamma^0}\right)^{-1}\right)^{-1}\le \Sigma_{\gamma^0}^{-1}$, we have
\begin{eqnarray*}
&&\frac{1+\textbf{y}^T (I_n-X_{\gamma^0} U_{\gamma^0}^{-1} X_{\gamma^0}^T )\textbf{y}}{1+\textbf{y}^T (I_n-P_{\gamma^0})\textbf{y}}\\
&=&1+\frac{\Delta}{1+\textbf{y}^T(I_n-P_{\gamma^0})\textbf{y}} \\
&\le&1+2\left(\frac{(\beta^0_{\gamma^0})^T\Sigma_{\gamma^0}^{-1}\beta^0_{\gamma^0}+\epsilon^TX_{\gamma^0}(X_{\gamma^0}^TX_{\gamma^0})^{-1} \Sigma_{\gamma^0}^{-1}(X_{\gamma^0}^TX_{\gamma^0})^{-1}X_{\gamma^0}^T\epsilon}{1+\textbf{y}^T(I_n-P_{\gamma^0})\textbf{y}}\right)\\
&\le&1+2\underline{\phi}_n^{-1}\left(\frac{\|\beta^0_{\gamma^0}\|_2^2 +\epsilon^TX_{\gamma^0}(X_{\gamma^0}^TX_{\gamma^0})^{-2} X_{\gamma^0}^T\epsilon} {1+\textbf{y}^T(I_n-P_{\gamma^0})\textbf{y}}\right).
\end{eqnarray*}
Since $\textbf{y}^T(I_n-P_{\gamma^0})\textbf{y}/n=\epsilon^T(I_n-P_{\gamma^0})\epsilon/n\rightarrow_p\sigma_0^2$, and $E\{\epsilon^TX_{\gamma^0}(X_{\gamma^0}^TX_{\gamma^0})^{-2} X_{\gamma^0}^T\epsilon\}\le \sigma_0^2s_n(n\varphi_{\min}(n))^{-1}$, we have
$\epsilon^T X_{\gamma^0}(X_{\gamma^0}^T X_{\gamma^0})^{-2} X_{\gamma^0}^T \epsilon=O_p\left(s_n(n\varphi_{\min}(n))^{-1}\right)$. Therefore,
by Assumptions \ref{assump:2} and \ref{assump:3}, and the fact that $k_n\ge s_n\psi_n^2$, we can show that
\begin{equation}\label{eq2}
\frac{1+\textbf{y}^T(I_n-X_{\gamma^0} U_{\gamma^0}^{-1}X_{\gamma^0}^T )\textbf{y}}{1+\textbf{y}^T(I_n-P_{\gamma^0})\textbf{y}} \le 1+\frac{2k_n}{n\underline{\phi}_n\sigma_0^2}(1+o_p(1)).
\end{equation}
Consequently, $0\le-T_4=O_p(1)$ follows from the condition that $k_n=O(\underline{\phi}_n)$ (Assumption \ref{assump:7}).

Next we approximate $T_2$ and $T_5$ in the following Lemmas \ref{lemma3} and \ref{lemma4}.

\begin{Lemma}\label{lemma3}
Under Assumption \ref{assump:8}, if $\gamma\in S_1$, then uniformly for $c_j$'s $\in[\underline{\phi}_n,\bar{\phi}_n]$,
$T_2\ge 2^{-1} (|\gamma|-s_n)\log(1+C_3n^{1-\delta}\underline{\phi}_n)$. Under
Assumption \ref{assump:2}, if $\gamma\in S_2$, then uniformly for $c_j$'s $\in[\underline{\phi}_n,\bar{\phi}_n]$,
$T_2\ge -2^{-1} s_n\log(1+C_2n\bar{\phi}_n)$, 
where $C_2$ and $C_3$ are constants given in Assumptions \ref{assump:2} and \ref{assump:8} respectively.
\end{Lemma}

{\bf Proof of Lemma \ref{lemma3}.}\,\, If $\gamma\in S_1$, it follows from the determinant formula for block matrices 
(Seber and Lee, 2003, page 468), and Assumption \ref{assump:8} that
\begin{eqnarray*}
\det(U_\gamma)&=&\det(U_{\gamma^0}) \det\left(\Sigma^{-1}_{\gamma\backslash\gamma^0}+X_{\gamma\backslash\gamma^0}^T
(I_n-X_{\gamma^0}U_{\gamma^0}^{-1}X_{\gamma^0}^T)X_{\gamma\backslash\gamma^0}\right)\\
&\ge&\det(U_{\gamma^0})
\det\left(\Sigma^{-1}_{\gamma\backslash\gamma^0}+X_{\gamma\backslash\gamma^0}^T
(I_n-P_{\gamma^0})X_{\gamma\backslash\gamma^0}\right)\\
&\ge&\det(U_{\gamma^0})\det\left(\Sigma_{\gamma\backslash\gamma^0}^{-1}+C_3n^{1-\delta}I_{|\gamma\backslash\gamma^0|}\right).
\end{eqnarray*}
Therefore,
\begin{eqnarray}\label{eq4}
\frac{\det(W_\gamma)}{\det(W_{\gamma^0})}&=&\frac{\det(\Sigma_\gamma)}{\det(\Sigma_{\gamma^0})}
\frac{\det(U_\gamma)}{\det(U_{\gamma^0})}\nonumber\\
&\ge&\det(\Sigma_{\gamma\backslash\gamma^0})\det\left(\Sigma_{\gamma\backslash\gamma^0}^{-1}+
C_3n^{1-\delta}I_{|\gamma\backslash\gamma^0|}\right)\nonumber\\
&=&\det\left(I_{|\gamma\backslash\gamma^0|}+C_3n^{1-\delta}\Sigma_{\gamma\backslash\gamma^0}\right)\nonumber\\
&\ge&\det\left((1+C_3n^{1-\delta}\underline{\phi}_n)I_{|\gamma\backslash\gamma^0|}\right) 
= (1+C_3n^{1-\delta}\underline{\phi}_n)^{|\gamma|-s_n},
\end{eqnarray}
which shows that $T_2\ge
2^{-1} (|\gamma|-s_n)\log(1+C_3n^{1-\delta}\underline{\phi}_n)$.
If $\gamma\in S_2$, note that $\det(W_{\gamma})\ge1$, and by Assumption \ref{assump:2}
\[
T_2\ge-\frac{1}{2}\log(\det(W_{\gamma^0}))\ge -\frac{1}{2}\log(\det(I_{s_n}+C_2n\Sigma_{\gamma^0}))\ge-2^{-1} s_n\log(1+C_2n\bar{\phi}_n),
\]
which completes the proof of Lemma \ref{lemma3}.
\qed

\begin{Lemma}\label{lemma4}
Let $\alpha_0>2$. If either Assumption \ref{assump:4} or \ref{assump:5} is satisfied, when $n$ is large, with large probability and uniformly for $\gamma\in S_1$, $T_5\ge-2^{-1} (|\gamma|-s_n)\alpha_0\log{p}$. If both Assumptions \ref{assump:2} and \ref{assump:4} are satisfied, there exists a constant $C'$ such that when $n$ is large, with large probability and uniformly for $\gamma\in S_2$, $T_5\ge 2^{-1}(n+\nu)\log\left(1+C' \psi_n^2\right)$.
\end{Lemma}

{\bf Proof of Lemma \ref{lemma4}.}\,\, We consider $\gamma\in S_1$ and $S_2$ separately. Notice that Assumption \ref{assump:4} implies that $p\log{p}=o(n\log(1+\psi_n^2))$, and therefore implies that $p\log{p}=o(n\psi_n^2)$. Let $v_\gamma=(I_n-P_\gamma)X_{\gamma^0\backslash\gamma}\beta^0_{\gamma^0\backslash\gamma}$. From Lemma \ref{lemma1} (a) and (c), there exists $C>0$ such that when $n$ is sufficiently large, with large probability, for any $\gamma\in S_2$,
\begin{eqnarray}\label{eq11}
\textbf{y}^T(I_n-P_\gamma)\textbf{y}&=&\|v_\gamma\|_2^2+2v_\gamma^T\epsilon+\epsilon^T(I_n-P_\gamma)\epsilon\nonumber\\
&\ge&\|v_\gamma\|_2^2-2C\sqrt{p}\|v_\gamma\|_2+\epsilon^T\epsilon-C|\gamma|\log{p}\nonumber\\
&\ge& \|v_\gamma\|_2^2\left(1-2C\frac{\sqrt{p}}{\|v_\gamma\|_2}-C\frac{p\log{p}}{\|v_\gamma\|_2^2}\right)+\epsilon^T\epsilon\nonumber\\
&\ge& \|v_\gamma\|_2^2\left(1-2C\sqrt{\frac{p}{n\varphi_{\min}(n)\psi_n^2}}-C\frac{p\log{p}}{n\varphi_{\min}(n)\psi_n^2}\right)+\epsilon^T\epsilon\nonumber\\
&=&\|v_\gamma\|_2^2(1+o(1))+\epsilon^T\epsilon\nonumber\\
&\ge&n\varphi_{\min}(n)\|\beta^0_{\gamma^0\backslash\gamma}\|_2^2(1+o(1))+\epsilon^T\epsilon\nonumber\\
&\ge& n\varphi_{\min}(n)\psi_n^2(1+o(1))+\epsilon^T\epsilon.
\end{eqnarray}
It is easy to see that Assumption \ref{assump:4} implies that $s_n=o(n)$, and therefore, $\epsilon^T(I_n-P_{\gamma^0})\epsilon=n\sigma_0^2 (1+o_p(1))$. Thus, by (\ref{eq11}), there exists a $C'$ such that for sufficiently large $n$, with large probability, uniformly for $\gamma\in S_2$,
\begin{equation}\label{eq1}
T_5\ge\frac{n+\nu}{2}\log\left(\frac{1+n\varphi_{\min}(n)\psi_n^2(1+o(1))
+\epsilon^T\epsilon}{1+\epsilon^T(I_n-P_{\gamma^0})\epsilon}\right)\ge\frac{n+\nu}{2}\log\left(1+C'\psi_n^2\right).
\end{equation}

On the other hand, by properties of projection matrices and Lemma \ref{lemma1} (b), when $n$ is sufficiently large,
 with large probability, we have uniformly for $\gamma\in S_1$,
\begin{eqnarray*}\label{eq12}
&&\frac{1+\textbf{y}^T(I_n-P_\gamma)\textbf{y}}{1+\textbf{y}^T(I_n-P_{\gamma^0})\textbf{y}}\\
&=&1-\frac{\textbf{y}^T(P_\gamma-P_{\gamma^0})\textbf{y}}{1+\textbf{y}^T(I_n-P_{\gamma^0})\textbf{y}}\nonumber\\
&=&1-\frac{(\beta_{\gamma^0}^0)^TX_{\gamma^0}^T(P_\gamma-P_{\gamma^0})X_{\gamma^0}\beta_{\gamma^0}
+2(\beta_{\gamma^0}^0)^TX_{\gamma^0}^T(P_\gamma-P_{\gamma^0})\epsilon+\epsilon^T(P_\gamma-P_{\gamma^0})\epsilon}
{1+\textbf{y}^T(I_n-P_{\gamma^0})\textbf{y}}\nonumber\\
&=&1-\frac{\epsilon^T(P_\gamma-P_{\gamma^0})\epsilon}
{1+\epsilon^T(I_n-P_{\gamma^0})\epsilon}\ge1-\frac{\alpha(|\gamma|-s_n)\log{p}}{n},
\end{eqnarray*}
where we have temporarily fixed an $\alpha$ such that $2<\alpha<\sqrt{2\alpha_0}$. 
It follows by the inequality that $\log(1-x)\ge -(\alpha/2)x$ when $x\in (0,1-2/\alpha)$, 
and by Assumption \ref{assump:4} or \ref{assump:5} (which both imply that $(|\gamma|-s_n)\log{p}/n$ 
approaches zero uniformly for $\gamma\in S_1$) that for sufficiently large $n$,
 with large probability and uniformly for $\gamma\in S_1$,
\begin{equation}\label{eq3}
T_5\ge\frac{n+\nu}{2}\log\left(1-\frac{\alpha(|\gamma|-s_n)\log{p}}{n}\right)\ge -2^{-1} (|\gamma|-s_n)\alpha_0\log{p},
\end{equation}
which completes the proof of Lemma \ref{lemma4}.
\qed

Now we are ready to finish the proof of Theorem \ref{thm1}. 
By (\ref{eq2}), Lemma \ref{lemma3}, Lemma \ref{lemma4}, 
Assumption \ref{assump:4}, and the fact that $p^{\alpha_0}=o(\rho_n)$ 
with $\rho_n\equiv n^{1-\delta}\underline{\phi}_n$, with large
probability, uniformly for $\gamma\in S_1$ and $c_1,\ldots,c_p\in[\underline{\phi}_n,\bar{\phi}_n]$,
\begin{eqnarray}\label{eq15}
p(\gamma|Z)/p(\gamma^0|Z)
&\le& \tilde{C}\,\,\exp\left(-2^{-1} (|\gamma|-s_n)\log((1+C_3\rho_n)/p^{\alpha_0})\right)\nonumber\\
&=&\tilde{C}
\left(\frac{1+C_3\rho_n}{p^{\alpha_0}}\right)^{-2^{-1} (|\gamma|-s_n)}\rightarrow0.
\end{eqnarray}
By Assumptions \ref{assump:4} and \ref{assump:6}, it can be verified that $s_n\log(1+C_2n\bar{\phi}_n)\ll\frac{n+\nu}{2}\log(1+C'\psi_n^2)$.
So, with large probability, uniformly for $\gamma\in S_2$ and $c_1,\ldots,c_p\in[\underline{\phi}_n,\bar{\phi}_n]$,
\begin{eqnarray}\label{eq16}
p(\gamma|Z)/p(\gamma^0|Z) 
&\le& \tilde{C}\,\,\exp\left(2^{-1} s_n\log(1+C_2n\bar{\phi}_n)-\frac{n+\nu}{2}\log(1+C'\psi_n^2)\right)\nonumber\\
&\le&\tilde{C}\,\,(1+C'\psi_n^2)^{-\frac{n+\nu}{4}}\rightarrow0,
\end{eqnarray}
where $\tilde{C}$ in (\ref{eq15}) and (\ref{eq16}) depends on the lower bounds of $T_1$ and $T_4$.
For the proof of PMC, we consider two cases. It is easy to see from (\ref{eq15}) that
\begin{eqnarray}\label{eq7}
\sum\limits_{\gamma\in
S_1}p(\gamma|Z)/p(\gamma^0|Z)&\le&\tilde{C}\sum\limits_{\gamma\in
S_1}\left(\frac{1+C_3\rho_n}{p^{\alpha_0}}\right)^{-2^{-1} (|\gamma|-s_n)}\nonumber\\
&=&\tilde{C}\sum\limits_{r=1}^{p-s_n}{p-s_n\choose r}
\left(\frac{1+C_3\rho_n}{p^{\alpha_0}}\right)^{-\frac{r}{2}}\nonumber\\
&=&\tilde{C}\left[\left(1+\left(\frac{1+C_3\rho_n}{p^{\alpha_0}}\right)^{-\frac{1}{2}}\right)^{p-s_n}-1\right]\rightarrow0\nonumber,
\end{eqnarray}
where the last limit result follows from the assumption that $p^{\alpha_0+2}=o(\rho_n)$.

Similarly, by (\ref{eq16}), and $p\log{n}=o(n\log(1+\psi_n^2))$ (which follows from Assumption \ref{assump:4}), we can show that
\begin{equation}\label{eq8}
\sum\limits_{\gamma\in S_2}p(\gamma|Z)/p(\gamma^0|Z)\le \tilde{C}2^p (1+C'\psi_n^2)^{-(n+\nu)/4}\rightarrow0.
\end{equation}
This completes the proof of Theorem \ref{thm1}.
\qed

{\bf Proof of Theorem \ref{thm3}.}\,\, The assumption that $\gamma^0$ is null implies that the model class $S_2$ is empty. 
Similar to the proof of Theorem \ref{thm1}, we need to approximate $T_1$ to $T_5$ in (\ref{eq5}). 
This is easier when the true model is null since $T_4=0$, and by Lemma
\ref{lemma3}, when $\gamma$ is nonnull, $T_2\ge 2^{-1}|\gamma|\log(1+C_3n^{1-\delta}\underline{\phi}_n)$. 
Since $T_1$ and $T_3$ are still bounded below, the proof is reduced to approximate $T_5$.
 By Lemma \ref{lemma4},  Assumption \ref{assump:5}, and that $s_n=0$, when $n$ is large, 
 with large probability and uniformly for $\gamma\in S_1$, $T_5\ge-2^{-1}|\gamma|\alpha_0\log{p}$. 
 Therefore, the remaining proofs can be finished by arguments similar to
(\ref{eq15}) and (\ref{eq7}).
\qed

{\bf Proof Theorem \ref{thm:g-prior:1}.}\,\, (i) Let $p(\gamma|Z,c)$ be the posterior probability of
$\gamma$ given $Z$ and $c$, as specified by (\ref{eq:post}). Applying Theorem \ref{thm1}, we have that in probability
\[
\sup\limits_{c\in[\underline{\phi},\bar{\phi}]}\max\limits_{\gamma\neq\gamma^0} p(\gamma|Z,c)/p(\gamma^0|Z,c)
\rightarrow 0.
\] 
Then the result follows from $p(\gamma|Z)=\int_{\underline{\phi}}^{\bar{\phi}} p(\gamma|Z,c) g(c)dc$, and
\[
\frac{\int_{\underline{\phi}}^{\bar{\phi}} p(\gamma|Z,c)g(c)dc}{\int_{\underline{\phi}}^{\bar{\phi}} p(\gamma^0|Z,c)g(c)dc}
\le \sup\limits_{c\in[\underline{\phi},\bar{\phi}]}\max\limits_{\gamma\neq\gamma^0} p(\gamma|Z,c)/p(\gamma^0|Z,c).
\]

(ii) Let $0<\underline{\phi}<\bar{\phi}$. By Theorem \ref{thm1},
$\inf\limits_{c\in[\underline{\phi},\bar{\phi}]} p(\gamma^0|Z,c)\rightarrow 1$ in probability.
Since
\[
p(\gamma^0|Z)=\int_{\underline{\phi}}^{\bar{\phi}} (p(\gamma^0|Z,c)-1) g(c)dc+\int_{\underline{\phi}}^{\bar{\phi}} g(c)dc
+\int_{[0,\infty)\backslash [\underline{\phi},\bar{\phi}]} p(\gamma^0|Z,c)g(c)dc,
\]
the result follows by fixing $\underline{\phi}$ and $\bar{\phi}$ so that $\int_{\underline{\phi}}^{\bar{\phi}} g(c)dc$ is close to 1,
and letting $n$ go to $\infty$.

{\bf Proof Theorem \ref{thm:g-prior:2}.} Proof is similar to those of Theorem \ref{thm:g-prior:1}.

\vspace{5mm}
\begin{center}
\textbf{Supplement Materials}
\end{center}

Supplements A--C are given in the authors' website: 
\begin{center}
http://www.stat.wisc.edu/$\sim$ shang/\vspace{5mm}
\end{center}
\textit{Supplement A:} Generalizations of Bayesian consistency to ultra-high dimensional settings.

\textit{Supplement B:} Proof of Corollaries \ref{cor:1} and \ref{cor:2}.

\textit{Supplement C:} Almost Sure Consistency of $p(\gamma^0|Z)$.

\vspace{5mm}
\noindent{\bf Acknowledgement}\,\, The authors wish to thank Professor 
Jun Shao for suggestions that helped to improve the present work, and
an anonymous reviewer who suggested the extension to $g$-priors.

\vspace{10mm}
\noindent{\large\bf References}
\begin{description}
\item Berger, J. O. and Pericchi, L. (1996). 
The intrinsic {B}ayes factor for model selection and prediction.
{\it J. Amer. Statist. Assoc.}
{\bf 91}, 109--122.

\item Berger, J. O., Ghosh, J. K. and Mukhopadhyay, N. (2003).
Approximations and consistency of Bayes factors as model dimension grows.
{\it J. Statist. Planning. Inference.}
{\bf 112}, 241--258.

\item B\"{u}hlmann, P., and Kalisch, M. and Maathuis, M. H. (2010).
Variable selection in high-dimensional linear models: partially faithful distributions and the
PC-simple algorithm. 
\textit{Biometrika} \textbf{97}, 261--278.

\item Casella, C., Gir{\'o}n, F. J., Mart\'{\i}nez, M.~L. and Moreno, E. (2009).
Consistency of {B}ayesian procedures for variable selection.
{\it Ann. Statist.}
{\bf 37}, 1207--1228.

\item Clyde, M. and George, E. I. (2000).
Flexible empirical {B}ayes estimation for wavelets.
{\it J. R. Stat. Soc. Ser. B.} {\bf 62}, 681--698.

\item Clyde, M., Parmigiani, G. and Vidakovic, B. (1998).
Multiple shrinkage and subset selection in wavelets.
{\it Biometrika} {\bf 85}, 391--401.

\item Durrett, R. (2005).
{\it Probability: Theorey and Examples}. 3rd Ed. Wadsworth-Brooks/Cole, Pacific Grove.

\item Fan, J. and Lv, J. (2008).
Sure independence screening for ultrahigh dimensional feature space.
{\it J. R. Stat. Soc. Ser. B.}
{\bf 70}, 849--911.

\item Fan, J. and Peng, H. (2004).
Nonconcave penalized likelihood with a diverging number of parameters.
{\it Ann. Statist.}
{\bf 32}, 928--961.

\item Fern{\'a}ndez, C.,  Ley, E. and Steel, M. F. J. (2001).
Benchmark priors for Bayesian model averaging.
{\it J. Econometrics}
{\bf 100}, 381--427.

\item George, E. and
McCulloch, R. (1993).
Variable selection via Gibbs sampling. 
\textit{J. Amer. Statist. Assoc.} \textbf{88}, 881--889.

\item Godsill, J.~S. and Rayner, P.~J.~W. (1998).
Robust reconstruction and analysis of autoregressive signals in impulsive noise using the Gibbs sampler. 
\textit{IEEE Trans. Speech Audio Process} \textbf{6}, 352--372.

\item Jeffreys, H. (1967).
{\it Theory of Probability}. 4th Ed. Oxford Univ. Press, Oxford.

\item Jiang, W. (2007).
Bayesian variable selection for high dimensional generalized linear models: Convergence rates of the fitted densities.
{\it Ann. Statist.}
{\bf 35}, 1487--1511.

\item Liang, F., Paulo, R., Molina, G., Clyde, M. and Berger, J. O. 
(2008).
Mixtures of $g$-priors for Bayesian variable selection.
{\it J. Amer. Statist. Assoc.}
{\bf 103}, 410--423.

\item Meinshausen, N. and B{\"u}hlmann, P. (2006).
High-dimensional graphs and variable selection with the lasso.
{\it Ann. Statist.}
{\bf 34}, 1436--1462.

\item Meinshausen, N. and Yu, B. (2009).
Lasso-type recovery of sparse representations for high-dimensional data.
{\it Ann. Statist.}
{\bf 37}, 246--270.

\item Moreno, E., Bertolino, F. and Racugno, W. (1998).
An intrinsic limiting procedure for model selection and hypotheses testing.
{\it J. Amer. Statist. Assoc.}
{\bf 93}, 1451--1460.

\item Moreno, E. and Gir{\'o}n, F. J. (2005).
Consistency of {B}ayes factors for intrinsic priors in normal linear models.
{\it C. R. Math. Acad. Sci. Paris}
{\bf 340}, 911--914.

\item Moreno, E., Gir\'{o}n, F. J. and Casella, G. (2010).
Consistency of objective Bayes factors as the model dimension grows. 
\textit{Ann. Statist.} \textbf{38}, 1937--1952.

\item Gir\'{o}n, F. J., Moreno, E., Casella, G. and
Mart\'{i}nez, M. L.
(2010). Consistency of objective Bayes factors for nonnested linear
models and increasing model dimension.
\textit{Revista de la Real Academia de Ciencias Exactas, Fisicas y Naturales. Serie A. Matematicas}
\textbf{104}, 57--67.

\item Shao, J. (2003).
{\it Mathematical Statistics}, 2nd Ed. Springer Texts in Statistics.
Springer, New York.

\item Seber, G.~A.~F. and Lee, A. J.
(2003).
{\it Linear Regression Analysis},
2nd Ed.
Wiley-Interscience [John Wiley \& Sons], Hoboken, NJ.

\item Smith, M. S. and Kohn, R. (1996).
Nonparametric regression using Bayesian variable selection. \textit{J. Econometrics} 
\textbf{75}, 317--344.

\item Wolfe, P. J., Godsill, S. J. and Ng, W.-J.
(2004).
Bayesian variable selection and regularization for time-frequency surface estimation.
{\it J. R. Stat. Soc. Ser. B.}
{\bf 66}, 575--589.

\item Zellner, A. (1971).
{\it An Introduction to {B}ayesian Inference in Econometrics}. Wiley, New York.

\item Zellner, A.(1978).
Jeffreys-Bayes posterior odds ratio and the Akaike information
criterion for discriminating between models.
{\it Econom. Lett.}{\bf 1}, 337--342.

\item Zellner, A. (1986). On assessing prior distributions and Bayesian regression analysis with g-prior
distributions. In \textit{Bayesian Inference and Decision Techniques: Essays in Honor of Bruno de
Finetti,} (eds. P. K. Goel and A. Zellner), 233--243. North-Holland/Elsevier.

\item Zhang, J., Clayton, M. K. and Townsend, P. (2010).
Functional concurrent linear regression model for spatial images. 
\textit{Journal of Agricultural, Biological and Environmental Statistics}, {\bf 16}, 105--130.

\item Zhang, C.-H. and Huang, J. (2008).
The sparsity and bias of the LASSO selection in high-dimensional linear regression.
{\it Ann. Statist.},
{\bf 36}, 1567--1594.

\end{description}
\end{document}